\documentclass[a4paper]{article}

\usepackage{graphics}

\usepackage{amsmath}
\usepackage{a4}
\usepackage{amssymb,epsfig}
\usepackage{color}

\usepackage[margin=1.25in]{geometry}

\numberwithin{equation}{section}

\newcommand{\la}{\left\langle}
\newcommand{\ra}{\right\rangle}

\newcommand{\lla}{\langle\!\langle}
\newcommand{\rra}{\rangle\!\rangle}
\newcommand{\p}{\partial}

\begin{document}

\title{Intrinsic Storage Valuation by Variational Analysis}
\author{Dmitry Lesnik}
\date{\today}

\maketitle

\begin{abstract}
The mathematical problem concerning intrinsic storage optimisation is formulated and solved by means of variational analysis. The solution, though obtained in implicit form, still sheds light on many important features of the optimal exercise strategy. It is shown how the solution depends on different constraint types including carry cost and cycle constraint.  Additionally, the relationship between intrinsic and stochastic solutions is investigated. In particular, we show that the optimal stochastic exercise decision is always close to the intrinsic one.
\end{abstract}

%\vspace{1cm}
%\begin{flushleft}
%D. Lesnik \\
%Gazprom Marketing \& Trading Limited, 20 Triton Street, \\ 
%London NW1 3BF, UK \\
%email: dmitry.lesnik@gmail.com\\
%\end{flushleft}
%

%\vspace{1cm}
%\begin{flushleft}
%\textbf{Keywords:} Storage option, Swing option, Optimal control,  Variational analysis, Pontryagin's principle
%\end{flushleft}
%
%\vspace{1cm}
%\begin{flushleft}
    %\textbf{Mathematics Subject Classification (2010):} 91G80,  49J30
%\end{flushleft}
%\begin{flushleft}
    %\textbf{JEL Classification:} G12 $\cdot$ G13
%\end{flushleft}

%\pagebreak

\vspace{1cm}

\section{Introduction}
\label{intro}

The problem of storage optimisation is driven by the necessity of the storage owners to minimise their expenses and maximise the potential profit which can be gained by operating the storage. There are plenty of real world storage examples: Gas storage, Oil storage, Hydro power plant, Coal stock, etc. Another popular product on the gas and electricity markets is a swing contract, which is mathematically equivalent to the storage problem. By following a clever strategy -- buying the underlying commodity cheap, storing it, and selling as the prices go up -- the storage owner can make a profit.  Knowing when, and how much, to buy or sell determines the exercise rule for the storage. Thus, there is a need for exercise strategy optimisation -- finding an exercise rule that maximises profit and minimises the risks.

An exact storage option pricing method must also take the stochastic nature of the prices into account. This is usually done by modelling prices with a stochastic process that, as far as possible, captures the value-generating statistical properties of the real world price movements. This price model is normally calibrated on the available market data before it being used for the actual pricing. We refer to the corresponding option value as stochastic value. Alternatively one may consider the storage option in the case when the prices are fixed (to the forward prices observed on the valuation date). The corresponding option value is named ``intrinsic''. The intrinsic value of the storage option is necessarily smaller than the stochastic value, since it does not take the time value into account.

One of the most popular ways of storage optimisation is the ``dynamic programming'' algorithm. It can be used to optimize both the intrinsic and stochastic problems. Stochastic methods approximate stochastic prices by modelling either on a tree or by Monte Carlo simulation. A very efficient numerical algorithm for valuation using Monte-Carlo pricing was suggested by Longstaff and Schwartz \cite{Schwartz2001}. A comprehensive description of the Least Squares Monte Carlo (LSMC) algorithm applied to the Storage Option is given in the paper by Alexander Boogert and Cyriel de Jong \cite{LSMC_Cyriel}. Dynamic programming offers the most generic approach, as it can deal with virtually any type of constraints. If the constraints are sufficiently simple then the intrinsic problem can be reformulated as a ``linear programming'' problem, which has much more efficient numerical algorithm, although usually the numerical efficiency is not an issue for the intrinsic problem.

The numerical solution of the storage option provides answers to the two most important questions: what is the optimal exercise decision at every time moment given the current state (prices and volume level of the storage), and what is the expected profit, provided we follow the optimal exercise strategy. Another important task is to find an optimal hedge strategy -- a portfolio of derivatives (futures, options or any other financial instruments) -- which would minimise the financial risks. Of course, market risk is only a feature of the stochastic problem, as the solution of the intrinsic problem is deterministic.  Non-market risks, such as operating risks of a physical storage, are not taken into account as they represent systemic risk, and so are in practice determined \emph{ex post facto}.

In this paper we develop an analytical approach to the intrinsic storage optimisation problem, and consider a link between intrinsic and stochastic solutions. The solution is obtained in implicit form, and hence can not be used for the option value calculation directly. However it highlights many interesting features of the problem, and covers a large number of applications. One of the most remarkable results is related to the optimal exercise decision. We show that the intrinsic, as well as stochastic, exercise decisions are based on a so-called trigger price, and most importantly -- the intrinsic and stochastic trigger prices coincide in the leading order. This implies that a correct stochastic exercise decision can be based on the intrinsic trigger price. Although a ``rolling intrinsic'' strategy is often used by traders, it is commonly believed to be sub-optimal. This work shows that the intrinsic exercise is in fact optimal, and we also demonstrate this with a numerical simulation.

Our analysis also shows how different constraints can be modelled in an efficient way. Among the most interesting examples is a cycle constraint, which is sometimes imposed on storage contracts. Treating a cycle constraint numerically is very computationally expensive task, as it requires an additional dimension in the state space. We show that the cycle constraint can be mimicked by a virtual injection/withdrawal costs, which makes the numerical solution much simpler to implement, and is by orders of magnitude faster.

A more detailed overview of the results is given in the section ''Discussion and numerical analysis``.

The paper is organised a follows. In Sec.~\ref{sec:ProblemFromulation} we give a formal mathematical formulation of the problem in terms of variational analysis. In Sec.~\ref{sec:Deterministic_problem} we present the solution to the intrinsic problem, and consider different special cases of constraints in Sec.~\ref{sec:other_constraint_types}. In Sec.~\ref{sec:StochTriggerPrice} we briefly consider certain aspects of the stochastic problem, and give a proof that the optimal stochastic exercise is bang-bang, and can be approximated by the intrinsic exercise. A reader who wishes to skip the detailed mathematical derivations can go straight to the discussion of the results in Sec.~\ref{sec:discussion}.

\section{Problem formulation}
\label{sec:ProblemFromulation}

The storage problem can be formulated as follows: The storage option holder is given a right to store some amount of an underlying asset (for the current discussion let it be natural gas) in a (virtual) storage facility. At each moment in time the option holder may ``do nothing'', inject gas into the storage, or release gas from the storage. Every time the gas is to be injected into the storage it must be bought on the market at the current spot price. Likewise, every time gas is released from the storage it is sold on the market. Since the market price of gas changes with time, this may lead to a non-trivial cash-flow

The injection and release of gas may have to satisfy some operational constraints (for instance maximum injection and release rates, maximum storage capacity, etc.) specified as boundary conditions. Every exercise profile (trajectory in the time-volume space) yields a different profit. In this sense, the profit becomes a functional on the exercise trajectory, mapping each possible profile of injection and release decisions to a corresponding profit. The aim of the storage option holder is to maximise the profit by choosing an optimal exercise strategy. The problem can thus be formulated in terms of variational analysis -- an optimal trajectory is the one delivering the maximum of the profit functional.

If the market prices are static, that is, the prices for all delivery days are known in advance and never change, then the problem is deterministic. The possible profit is bound from above and from below, and thus there exists a trajectory such that no other trajectory yields higher profit. The maximal profit of the deterministic problem is called ``intrinsic value''. We shall now investigate the deterministic problem by means of variational analysis.

We consider the storage problem in a continuous time approximation. Let $T_e$ be the total time of the storage contract, such that $t\in[0, T_e]$. 

Let $F(t)$ be the market price of gas for delivery at the time $t$. The curve $F(t)$ is called the forward price curve, since it is observed on the market at the time $t=0$ and contains the information about the prices of gas with delivery in the future. By the definition of the static problem the forward curve never changes, and hence there is no difference between observation and delivery time\footnote{In reality the forward curve changes during the trading, and thus it is a function of both observation time $t$ and delivery time $T$: $F(t,T)$.}.

Let $q(t)$ be the amount of gas contained in the storage at time $t$. We refer to the volume injected into (or released from) the storage, per unit time, as the ``exercise''. The curve $\dot q(t) = dq/dt$ defines the exercise trajectory. The initial and terminal conditions are
\begin{align}
     &q(0)= Q_{start}  \qquad \qquad   q(T_e)= Q_{end}
\end{align}
We consider the case of the free terminal condition, where $q(T_e)$ is not fixed, in Sec.~\ref{sec:free_vs_fixed_terminal_condition}.

The cash flow resulting from trading the gas according to the exercise strategy $\dot q(t)$ over the time interval $dt$ is given by
\begin{eqnarray}
    \label{eq:1.1}
    - \dot q(t)\,F(t)\, dt
\end{eqnarray}
An additional cash flow can occur from the costs of injecting or releasing gas (operating costs). Let us designate~$\gamma(\dot q)$ as the operating costs per unit time (we assume here that the operating cost is only a function of the exercise $\dot q$). The terminal profit or loss is given by the cumulative cash flow over the storage lifetime. We thus introduce the target functional
\begin{eqnarray}
    \label{eq:unperturbed_action1}
    &&S_0[q(t)] = -\int_0^{T_e} \Big[ \dot q(t)\,F(t) + \gamma(\dot q(t)) \Big]\,dt\\
    \label{eq:unperturbed_action2}
    &&q(0)= Q_{start}\,,  \qquad    q(T_e)= Q_{end}
\end{eqnarray}
The value of this functional on the trajectory $q(t)$ gives the storage value conditional on that trajectory occurring. Here we have used the subscript $0$ to indicate the unmodified integral. In the next section we will introduce a modified integral $S[q(t)]$, which includes additional terms intended to enforce the operational constraints. Looking ahead, we notice that on any trajectory $q(t)$ allowed by the constraints the values of the modified and unmodified functionals coincide.

Using ``physics'' terminology we introduce the (unmodified) Lagrangian $L_0$ corresponding to the action functional $S_0$ as
\begin{align}
    \label{eq:unperturbed_lagrangian}
    L_0 = - \dot q(t)\,F(t) - \gamma(\dot q(t))
\end{align}

As mentioned above, the target functional is bounded both above and below, and hence there must exist trajectories on which the functional achieves its maximum and 
minimum values. We now make use of variational analysis to search for the extremal trajectory, i.e. the trajectory on which the first variation of the target functional vanishes\footnote{We do not conduct a thorough analysis of existence and uniqueness of the solution of the variational problem. However, we point out that the class of functions $q(t)$ on which the target functional is defined and achieves a maximum is rather broad. In particular, this class includes all continuous locally integrable functions. In some cases the solution space can even be extended to discontinuous or singular functions. We assume that all functions representing problem parameters -- forward curve, boundary conditions, etc. -- are sufficiently well behaved such that all mathematical expressions are well defined. If this is not the case, then we can apply the procedure of regularisation. We also utilise the concept of convergence of functions, which we will always understand to be weak convergence, i.e. $\psi_k\to\psi$ if $S[\psi_k]\to S[\psi]$.}.

The types of operational constraints may differ depending on the different storage option types. Below we consider some typical constraints found in practice, which can be classified as \emph{local} constraints. A local property implies that a constraint at time $t$ can be expressed completely in terms of state variables and their derivatives $q(t),\dot q(t), F(t),\dot F(t), ... $ at each time $t$.

We consider the following two operational constraints:
\begin{enumerate}
    \item The volume $q(t)$ must remain within the interval
    \begin{eqnarray}
        q(t) \in [Q_{min}(t),Q_{max}(t)]\
    \end{eqnarray}
    where the boundaries $Q_{min}(t),Q_{max}(t)$ can be time dependent.
    \item The injection/release rate $\dot q(t)$ is bounded by
    \begin{eqnarray}
        \dot q(t) \in [r_{min}, r_{max}]
    \end{eqnarray}
    Generally the maximal injection/release rates may depend on both the time and volume: $r = r(t,q(t))$. In the following analysis we will only consider constant injection/release rates. The case of volume dependent rates will be considered in Sec.~\ref{sec:volume_dependent_rates}.
\end{enumerate}

One of the typical examples of a \emph{non-local} constraint is the so called \emph{cycle} constraint. It can be formulated as follows: One introduces an \emph{intake} cycle variable as
\begin{align}
    \label{eq:cycle_variable}
    c(t) = \int_0^{t} \dot q(\tau)\, \theta( \dot q(\tau) )\,d\tau
    \qquad\text{where} \quad
    \theta(x) = \left\{
        \begin{array}{ll}
            1\quad & x\geq 0 \\
            0      & x < 0
        \end{array}
    \right.
\end{align}
which corresponds to the total injected volume up to time $t$. The cycle constraint requires that the terminal value $c(T_e)$ does not exceed a certain threshold
\begin{align}
    \label{eq:cycle_constrain}
    c(T_e) \leq c_{max}
\end{align}
Similarly one can introduce a \emph{release} cycle constraint, bounding the total released volume. We will consider the cycle constraint below in section~\ref{sec:cycle_constraint}.

\subsection{Penalty functions}
\label{sec:PenaltyFunctions}

To restrict $q(t)$ from going beyond the range $[Q_{min},Q_{max}]$, we can introduce a parametrised penalty function $-\phi[q(t),N_\phi]$ and add it to the Lagrangian. The penalty function can be any smooth convex function, which in the limit $N_\phi\to \infty$ approaches zero within $[Q_{min},Q_{max}]$ and positive infinity otherwise. A particular example of this function could be
\begin{align*}
    & \phi = \Big[a(q(t)-b)\Big]^{2N_\phi} \\ 
    & \text{where} \quad
    a = \frac{2-\varepsilon_N}{Q_{max}-Q_{min}} \qquad
    b = \frac{Q_{max} + Q_{min}}{2}\qquad 
    \varepsilon_N = \frac{1}{N_\phi}
\end{align*}
The parametrised representation has a purpose of \emph{regularisation}, as the variational problem satisfies the strict requirements of classical variational analysis for finite $N_\phi$. The transition to the limit $N_\phi\to\infty$ should formally take place only after the finite $N_\phi$ solution is found.

A similar regularised penalty function $-\psi[\dot q(t),N_\psi]$ can be introduced to restrict~$\dot q(t)$ from exiting the interval $[r_{min}, r_{max}]$. In the limit $N_\psi\to\infty$ it approaches zero within $[r_{min}, r_{max}]$ and positive infinity otherwise. The modified Lagrangian becomes
\begin{align}
    \label{eq:perturbed_lagrangian}
    L = L_0 - \phi(q) - \psi(\dot q) = - \Big[ \dot q F + \phi(q) + \psi(\dot q) + \gamma(\dot q) \Big]
\end{align}

In the next section we find a formal solution to the variational problem, and show that the solutions of modified and unmodified constrained problems coincide.

\section{Solution of deterministic problem}
\label{sec:Deterministic_problem}

\subsection{Euler-Lagrange equation}

We know from variational analysis that the extremal trajectory of the integral $\int L(q, \dot q)\,dt$ satisfies the Euler-Lagrange equation (see e.g. \cite{LandauMechanics})
\begin{align}
    \frac{\p L(q, \dot q)}{\p q} = \frac{d}{dt}\frac{\p L(q,\dot q)}{\p \dot q}
\end{align}
which is valid for both free and fixed terminal conditions. For the modified Lagrangian we obtain the following Euler-Lagrange equation
\begin{align}
    \label{eq:EL}
    \phi'(q) = \frac{d}{dt}\Big[F + \psi'(\dot q) + \gamma'(\dot q) \Big]
\end{align}

There are two types of solution for the latter equation. One is obtained on the interval where the trajectory remains strictly within the boundaries
    \[ q(t) \in (Q_{min}, Q_{max} )\]
The second solution type is where the trajectory lies on the boundary
    \[ q(t) \equiv Q_{min}(t) \quad \text{ or } \quad q(t) \equiv Q_{max}(t)\]
The general solution consists of pieces of both the solutions types. Let us consider each of these types separately.

\subsection{Solution between the boundaries}

First we consider the solution in the interval where the trajectory does not touch the boundary. In the limit $N_\phi\to\infty$ the l.h.s. of Eq.~(\ref{eq:EL}) vanishes. Thus we have
\begin{align}
    \label{eq:within_boundaries}
    F(t) + \psi'(\dot q) + \gamma'(\dot q) = C
\end{align}
where $C$ is \textit{constant} for the whole period of time, where the solution does not touch the boundary.

Let us consider one particular example of the operating cost function. It is typical for a gas storage facility that the operating cost is proportional to the amount of the gas released or injected into the storage. Again we can parametrise it with $N_\gamma$ so that for finite $N_\gamma$ the function $\gamma(\dot q)$ is smooth, and in the limit $N_\gamma\to \infty$ it is piecewise linear
\begin{align}
    \gamma(\dot q) = \left\{
        \begin{array}{cc}
            \gamma_{inj} \, \dot q \qquad& \text{ for } \dot q > 0\,;\\
            -\gamma_{rel} \, \dot q \qquad& \text{ for } \dot q < 0\,;
        \end{array}
    \right.
\end{align}
where $\gamma_{inj} > 0$ and $\gamma_{rel}>0$.

The solution of Eq.~(\ref{eq:within_boundaries}) is now straightforward to find graphically (see Fig.~\ref{fig:1}). In case $r_{min}<0,\, r_{max}>0$ we obtain:
\begin{align}
    \label{eq:solution}
    \dot q(t) = \left\{
        \begin{array}{ll}
            r_{min}(t) \qquad &  F(t) > C+\gamma_{rel}               \\
            0                 &  C-\gamma_{inj} < F(t) < C+\gamma_{rel}  \\
            r_{max}(t) \qquad &  F(t) < C - \gamma_{inj}            
        \end{array}
    \right.
\end{align}
We see that the solution has three regimes, and the value $C$ can be interpreted as a trigger price. If the price $F(t)$ is within the ``dead zone'' $[C-\gamma_{inj},\, C+\gamma_{rel}]$, the optimal exercise is ``do nothing'': $\dot q = 0$. If the price is below the dead zone, the gas is injected at maximal rate: $\dot q(t) = r_{max}(t)$. If the price is above the dead zone, the gas is released at maximum rate: $\dot q(t) = r_{min}(t)$. This result is a generalisation of the well-known Pontryagin's maximum principle.

\setlength{\unitlength}{1mm}
\begin{figure}[htb!]
    \begin{center}
        \begin{picture}(90,60)
            \put(11,15){$r_{min}$}
            \put(68,15){$r_{max}$}
            \put(43.5,12){$-\gamma_{rel}$}
            \put(35,26){$\gamma_{inj}$}
            \put(10,37){$C - F(t)$}
            \put(83,15){$\dot q$}
            \put(71,48){$\psi'(\dot q) + \gamma'(\dot q)$}
%             \linethickness{1pt}
            \thicklines
            \put(2,20){\vector(1,0){83}}
            \put(42,2){\vector(0,1){53}}
            %\thinlines
            \put(12,19){\line(0,1){2}}
            \put(69,19){\line(0,1){2}}
            %\thicklines
            \color{red}
            \put(12,2){\line(0,1){11}}
            \put(12,13){\line(1,0){30}}
            \put(42,13){\line(0,1){14}}
            \put(42,27){\line(1,0){27}}
            \put(69,27){\line(0,1){25}}
            \color{blue}
            %\thicklines
            %\thinlines
            \put(5,35){\line(1,0){80}}
        \end{picture}
        \caption{\small Graphical illustration of Eq.~(\ref{eq:within_boundaries}).}
        \label{fig:1}
    \end{center}
\end{figure}
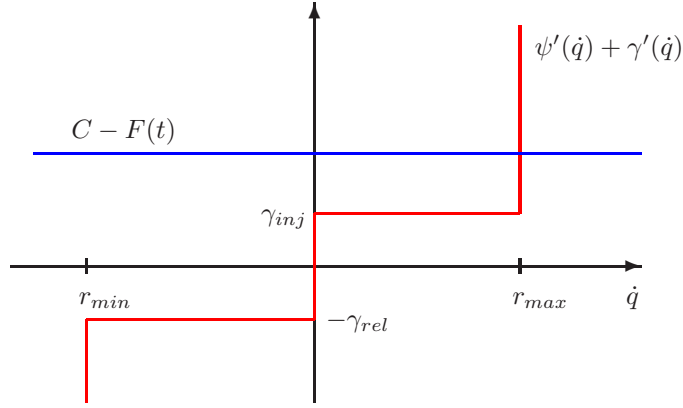

In particular, without operating costs the solution has only two regimes -- injection and release (this feature is usually called in literature ``bang-bang''). The transition between regimes occurs every time the forward curve crosses the trigger level $C$.

The constant $C$ has to be chosen in such a way that the solution $q(t)$ satisfies the initial and terminal conditions
\begin{align}
    \label{eq:init_term_cond}
    q(0) = Q_{start} \qquad \qquad q(T_e) = Q_{end} = Q_{start} + \int_0^{T_e} \dot q \, dt
\end{align}
and thus can not be represented analytically in the general case.

\subsection{Solution on the boundary}

There is not much we need to say about this part of the solution. On the boundary we obviously have
\begin{align*}
    & q(t) = Q_{max}(t)\,\qquad\text{on the upper boundary} \\
    & q(t) = Q_{min}(t)\,\qquad\text{on the lower boundary} 
\end{align*}
There are two different points of view on the trajectory lying on the boundary. We may consider our problem as a general variational problem with a spatial boundary condition. A solution, that maximises the functional, can consist of parts going along the boundary and \emph{extremals} -- pieces of trajectory lying between boundaries and satisfying the extremal condition (i.e. having vanishing first variation). The Euler-Lagrange equation does not make sense on the boundary, since the trajectory can not be varied there.

On the other hand, if we consider a regularised problem, the Euler-Lagrange equation does make sense on the boundary. We can derive the special case of the equation on the boundary. If the boundary $Q_{max}(t)$ (or $Q_{min}(t)$) is a sufficiently slow function of time, then the function $\psi'(\dot q)$ vanishes, and the function $\phi'(q)$ takes a finite value (positive if $q = Q_{max}$ or negative if $q = Q_{min}$). Eq.~(\ref{eq:EL}) on the boundary becomes
\begin{align}
    \label{eq:on_the_boundary}
    \phi'(q) = \frac{d}{dt} \Big[F(t)+\gamma'(\dot q)\Big]
\end{align}
From this follows 
\begin{align}
    & \frac{d}{dt} \Big[F(t)+\gamma'(\dot Q_{max})\Big]\geq 0\qquad
    \text{on the upper boundary;}\\
    & \frac{d}{dt} \Big[F(t)+\gamma'(\dot Q_{min})\Big]\leq 0\qquad
    \text{on the lower boundary;}
\end{align}
In particular for $Q_{max} = const$ and $Q_{min} = const$ we obtain
\begin{align}
    & \dot F(t) \geq 0\qquad \text{on the upper boundary;}\\
    & \dot F(t) \leq 0\qquad \text{on the lower boundary;}
\end{align}

\subsection{Approaching the boundary}

An important statement can be made about the connection point between extremals and the boundary. Let us first consider a simplified problem (a generalisation is straightforward), for which we suppose a constant boundary condition, and zero operating cost:
    \[ Q_{max} = const \qquad Q_{min} = const \qquad  \gamma(\dot q)\equiv 0\]
Let $q(t)$ be the optimal trajectory and $t^*$ be the time when the trajectory touches the boundary, such that, for $t<t^*$ the trajectory satisfies the Euler-Lagrange equation and for $t>t^*$ the trajectory follows the boundary, e.g. the upper boundary: $q(t>t^*) = Q_{max}$. 

Let us now consider a supplementary problem for the time horizon $t\in [0, t^*]$ with a terminal condition $q(t^*)=Q_{max}$. We designate $S_{sup}$ the value of the supplementary target function. Obviously, the optimal trajectory of the supplementary problem coincides on $[0,t^*]$ with the original optimal trajectory
\[q_{sup}(t) = q(t) \qquad  t\in [0,t^*] \]
Now consider a small variation of the boundary touch point
\[
    \tilde t^* = t^* + \delta t^*
\]
For the supplementary problem we will obtain a new optimal solution $\tilde q_{sup}(t)$.  As will be shown in Sec~\ref{sec:time_derivative}, the change of the value of the supplementary target function can be obtained from the time derivative of the boundary problem at the point $t = t^*$, which is given by
\begin{align}
    \label{eq:time_derivative_on_boundary}
    &\frac{\p S_{sup}}{\p t} = \dot q_{sup}\,\left( C - F(t^*)\right)\\
    \label{eq:time_derivative_on_boundary2}
    &\delta S_{sup} = \frac{\p S_{sup}}{\p t}\,\delta t^*
\end{align}
where $\dot q_{sup} = \dot q_{sup}(t^*-0)=r_{max}$ is the optimal exercise of the part of the trajectory just before the boundary touch.

Regarding the original problem, we can now ask how the value of the target function has changed on the new trajectory $\tilde q(t)$, which consists of $\tilde q(t)$ for $t < \tilde t^*$ and is unmodified for $t>\max( t^* , \tilde t^*)$ 
\[
    \tilde q(t) = \left\{
        \begin{array}{ll}
            \tilde q_{sup}(t)\,, \qquad &  t < \tilde t^*     \\
            q(t)                      &  t > \max( t^* , \tilde t^*)     \\
            Q_{max}                   &  \tilde t^* < t < t^* \quad \text{if}\quad \tilde t^* < t^*
        \end{array}
    \right.
\]
For the original variational problem the point $t^*$ is an internal point, and thus the new trajectory is a small variation of the optimal trajectory. Thus the original target functional remains unmodified:
\[
    \delta S = 0
\]
The change of the value of the original target function is only due to the change in the extremal part of the trajectory for $t<t^*$, since the part of the trajectory on the boundary does not contribute to the functional value. We thus conclude that the variation of the supplementary target function should vanish
\[
    \delta S_{sup} = 0
\]
Comparing this with Eqs.~(\ref{eq:time_derivative_on_boundary}) and~(\ref{eq:time_derivative_on_boundary2}) we can conclude that, at the time of boundary touch,
\begin{align}
    F(t^*) = C
\end{align}

In the general case (which includes time dependent boundary conditions and operating costs) we obtain the following conditions for the boundary touch points
\begin{enumerate}
    \item The trajectory touches the lower boundary from the left (Fig.~\ref{fig:boundary}.a):
        \begin{align}
            \label{eq:bc1}
            F(t^*) = C + \gamma_{rel}\,, \qquad \dot F(t^*) \leq 0
        \end{align}
    \item The trajectory touches the upper boundary from the left (Fig.~\ref{fig:boundary}.b):
        \begin{align}
            \label{eq:bc2}
            F(t^*) = C - \gamma_{inj}\,, \qquad \dot F(t^*) \geq 0
        \end{align}
    \item The trajectory touches the lower boundary from the right (Fig.~\ref{fig:boundary}.c):
        \begin{align}
            \label{eq:bc3}
            F(t^*) = C - \gamma_{inj}\,, \qquad \dot F(t^*) \leq 0
        \end{align}
    \item The trajectory touches the upper boundary from the right (Fig.~\ref{fig:boundary}.d):
        \begin{align}
            \label{eq:bc4}
            F(t^*) = C + \gamma_{rel}\,, \qquad \dot F(t^*) \geq 0
        \end{align}
\end{enumerate}
%\setlength{\unitlength}{1cm}
%\begin{figure}[htb!]
    %\begin{center}
        %\begin{picture}(8,2)
            %\thicklines
            %\put(-0.5,0.5){\vector(0,1){2.2}}
            %\put(-0.6,0.8){\line(1,0){9.5}}
            %\put(-0.6,1.8){\line(1,0){9.5}}
            %\thicklines
            %\color{red}
            %\put(-0.2,1.11){\line(2,-1){0.6}}
            %\put(0.4,0.81){\line(2,0){0.8}}
            %\put(2.3,1.49){\line(2,1){0.6}}
            %\put(2.9,1.79){\line(2,0){0.8}}
            %\put(4.6,0.81){\line(2,0){0.8}}
            %\put(5.4,0.81){\line(2,1){0.6}}
            %\put(6.9,1.79){\line(2,0){0.8}}
            %\put(7.7,1.79){\line(2,-1){0.6}}
            %\color{black}
            %\put(-1.2,2.5){$q(t)$}
            %\put(-1.5,1.7){$Q_{max}$}
            %\put(-1.5,0.7){$Q_{min}$}
            %\put(0.1,1.0){$a$}
            %\put(2.5,1.3){$b$}
            %\put(5.4,1.0){$c$}
            %\put(7.9,1.3){$d$}
        %\end{picture}
        %\vspace{-10pt}
        %\caption{\small Boundary touch cases: a) lower boundary from left; b) upper boundary from left; c) lower boundary from right; d) upper boundary from right.}
        %\label{fig:boundary}
    %\end{center}
%\end{figure}
\setlength{\unitlength}{1cm}
\begin{figure}[htb!]
    \begin{center}
        \begin{picture}(8,4)
            \put(-2,-0.2){\includegraphics[width=0.8\textwidth]{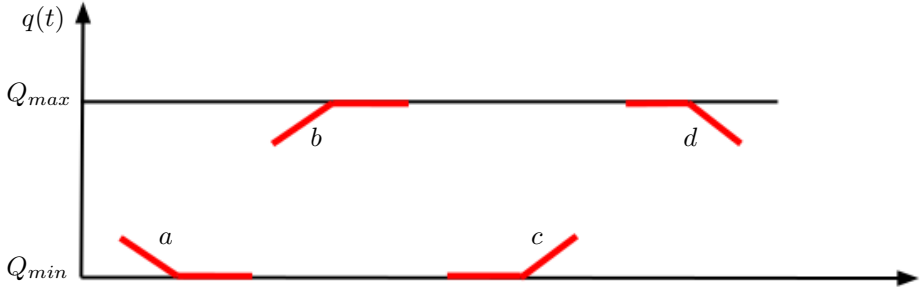}}
            \color{black}
            \put(-2.4,3.7){$q(t)$}
            \put(-2.6,2.7){$Q_{max}$}
            \put(-2.6,0.4){$Q_{min}$}
            \put(-0.6,0.8){$a$}
            \put(1.4,2.1){$b$}
            \put(4.3,0.8){$c$}
            \put(6.3,2.1){$d$}
        \end{picture}
        \caption{\small Boundary touch cases: a) lower boundary from left; b) upper boundary from left; c) lower boundary from right; d) upper boundary from right.}
        \label{fig:boundary}
    \end{center}
\end{figure}
These conditions must be satisfied for the time interval $t\in(0,T_e)$. The extremal trajectory at the times $t=0$ and $t=T_e$ does not have to satisfy these conditions as they are not internal points of the trajectory.

\subsection{Conclusion}

Eq.~(\ref{eq:solution}) together with the boundary conditions (\ref{eq:bc1}-\ref{eq:bc4}) and a condition on the start and end volume~(\ref{eq:init_term_cond}) provides an implicit solution for the problem. We summarise some properties of the solution:
\begin{enumerate}
    \item The optimal exercise strategy between boundaries is bang-bang, i.e. there are only 3 possible optimal actions at each time moment: maximum injection, maximum release or ``do nothing''. 
    \item Generally the optimal trajectory consists of a number of sections -- those following the boundary and interior sections satisfying Euler-Lagrange equation.
    \item For each piece of trajectory, separated from others by a boundary touch, there is a trigger price $C$ and a ``dead zone'' $[C-\gamma_{inj},\, C+\gamma_{rel}]$. The volume in the storage is increasing at maximum rate if the price is below the zone, decreasing at maximum rate if the price is above the zone, and is constant if the price is in the zone. The width of the zone is defined by the operating costs. Parts of the trajectory separated by a boundary touch can have different trigger levels.
    \item If the trajectory touches the boundary in the time interval $t\in(0,T_e)$, it must satisfy one of the boundary conditions~(\ref{eq:bc1}-\ref{eq:bc4}).
    \item Generally a trajectory satisfying all the previous conditions is not unique. Any solution must obey all these conditions, but not any function obeying all conditions is an optimal solution.
\end{enumerate}

The solution found in the previous section maximises the modified target functional $S[q(t)]$. Let $\bar q(t)$ be the obtained solution. The optimal path $\bar q(t)$ can be substituted to the unmodified functional $S_0[\bar q(t)]$, on which we expect to obtain the value of the option. There are two questions we need to address: 1) if the optimal trajectory $\bar q(t)$ respects the constraints, and 2) if the trajectory $\bar q(t)$, which is optimal for the modified functional, is also optimal for the unmodified one. For both these questions the answer is affirmative. In Appendix~\ref{sec:modified_unmodified_functionals} we show why.

From this point onwards we drop the subscript ``0'' for the action integral, unless we want to emphasise the difference between modified and unmodified target functions. We also drop the ``bar'' for the optimal trajectory $\bar q(t)$ everywhere where it does not lead to ambiguity.

\subsection{Dependency on initial and terminal boundary conditions}
\label{sec:3}

The optimal path $q(t)$ and the target functional $S[q(t)]$ are calculated for fixed terminal point $\{t=T_e,q=Q_{end}\}$. If this point is slightly shifted in either direction -- time or volume -- the optimal path becomes different, and the target function as well. In this sense the target function can be viewed as a function of the terminal point $S = S(q,t)$. In this section we find the derivatives of this function.

From the standard variational analysis we have:
\begin{align}
    \label{eq:term_derivative_space}
    & \frac{\partial S}{\partial q} = \frac{\partial L}{\partial \dot q} \\
    \label{eq:term_derivative_time}
    & \frac{\partial S}{\partial t} = L - \dot q \frac{\partial L}{\partial \dot q}
\end{align}
We apply these equations to the modified objective function $S[q(t)]$. The obtained formulas will also be valid for the unmodified objective function $S_0[q(t)]$.

\subsubsection{Volume derivative}
\label{sec:spacial_derivative}

Using the definition (\ref{eq:perturbed_lagrangian}) of the modified Lagrangian and Eq.~(\ref{eq:term_derivative_space}) we find
\[
    \frac{\partial S}{\partial q} = - \Big[ F + \psi'(\dot q) + \gamma'(\dot q) \Big]
\]
The spacial derivative makes sense only within the boundaries. Making use of Eq.~(\ref{eq:within_boundaries}) and replacing $S$ with $S_0$ we finally obtain
\begin{align}
    \label{eq:spacial_derivative}
    \frac{\partial S_0}{\partial q_{end}}  = - C
\end{align}
where $C$ is the trigger level of the last part of the trajectory. 

The derivative with respect to the initial condition $q_{start}$ is given by inverting the sign,
\begin{align}
    \frac{\partial S_0}{\partial q_{start}} = C
\end{align}
where $C$ is the trigger level of the initial part of the trajectory.

\subsubsection{Time derivative}
\label{sec:time_derivative}

From Eqs.~(\ref{eq:perturbed_lagrangian}) and~(\ref{eq:term_derivative_time}) we find
\begin{align}
    \frac{\partial S}{\partial t}& =
    \dot q\Big[\psi'(\dot q) + \gamma'(\dot q)\Big] -
        \phi(q) - \psi(\dot q) - \gamma(\dot q) = \nonumber \\
    & = \dot q\Big[\psi'(\dot q) + \gamma'(\dot q)\Big] - \gamma(\dot q)
\end{align}
\textbf{Within the boundaries} we make use of the relation~(\ref{eq:within_boundaries}). Replacing $S$ with $S_0$ we obtain
\begin{align}
    \frac{\partial S_0}{\partial t_{end}} = \dot q\, \Big[C - F(t_{end}) \Big] - \gamma(\dot q)
\end{align}
where $\dot q = \dot q(t_{end}-0)$ is the optimal exercise on the final part of the trajectory, obtained for a fixed terminal condition. In particular, one can see that  $\partial S_0/\partial t_{end} \geq 0$ (provided that $r_{min}<0, r_{max}>0$).

\textbf{On the boundary} in the special case $\dot Q_{max} = \dot Q_{min} = 0$ one obtains $\dot q= 0, \gamma(\dot q)= 0$ and
\begin{align}
    \frac{\partial S_0}{\partial t_{end}} = 0
\end{align}
It is worth noting that the derivative with respect to the terminal (or initial) boundary condition is universal: it is independent on whether the optimal trajectory touches the boundary or not. To calculate the derivatives, one just needs to know the solution $q(t)$ and the trigger price $C$ at the end (or beginning) of the trajectory. This also remains valid with imposed cycle constraint (Sec.~\ref{sec:cycle_constraint}).

\section{Other constraint types}
\label{sec:other_constraint_types}

So far we have only considered the optimisation problem with the simplest configuration. Often, realistic storage and swing contracts have additional constraints imposed, driven by operational or financial reasons. In this section we consider some typical additional constraints for storage and swing contracts, and discuss what influence they have on the solution.

\subsection{Carry cost}

Some storage contracts may include a ``carry cost'', which is the cost of keeping the underlying commodity in storage. For instance, an oil storage is subject to a heating cost (fuel oil in a storage has to be kept warm at approximately $60^oC$).

The carry cost is specified as a time-dependent price $\gamma_c(t)$ per unit time per unit volume. Thus, the (unmodified) target functional becomes
\begin{eqnarray}
    \label{eq:action_carry_cost}
    S_0 = -\int_0^{T_e} \Big[ \dot q\,F(t) + \gamma(\dot q) + \gamma_c(t)\,q(t) \Big]\,dt
\end{eqnarray}
Substituting the new Lagrangian into the Lagrange-Euler equation we obtain between boundaries,
\begin{align}
    \label{eq:within_boundaries_carry_cost}
    F(t) + \psi'(\dot q) + \gamma'(\dot q) = C(t)
\end{align}
where
\begin{align}
    \label{eq:trigger_price_carry_cost}
    \frac{dC(t)}{dt} = \gamma_c(t)
\end{align}
Thus, the solution of the problem with the carry cost is similar to that without, but has one difference: the trigger price is not a constant, rather, it is a growing function of time, satisfying Eq.~(\ref{eq:trigger_price_carry_cost}). 

It is easy to show that the carry cost can be replicated by an additional discount factor. Indeed, the trigger times occur when the forward curve crosses the trigger level $F(t) = C(t)$. Using the solution for Eq.~(\ref{eq:trigger_price_carry_cost}) in the form $C(t) = C_0 + \int \gamma_c(u) \, du$ we find the condition for the trigger time
\[
    F(t) - \int_0^t \gamma_c(u)\,du = C_0
\]
Thus the solution for the problem with carry cost is equivalent to the solution with modified forward curve $\tilde F(t) = D(t)\,F(t)$, where the discount factor $D(t)$ is given by
\begin{align}
    D(t) = 1 - \frac{1}{F(t)}\int_0^t \gamma_c(u)\,du \approx
    \exp\left(- \frac{1}{F(t)}\int_0^t \gamma_c(u)\,du\right)
\end{align}

\subsection{Solution with free terminal condition}
\label{sec:free_vs_fixed_terminal_condition}

Often the storage problem can be defined with a free terminal condition. For instance, one may be allowed to leave an arbitrary amount of the underlying in the storage, and for any remaining volume one gets a terminal pay-off, equivalent to selling that volume at an effective price $F_e$, referred to as the \emph{residual unit price}. We define the (unmodified) target functional $\tilde S$ as
\begin{align}
    \label{eq:free_term_cond}
    \tilde S[q(t)] = -\int_0^{T_e} \Big[ \dot q(t)\,F(t) + \gamma(\dot q)\Big]\,dt +
    q_{end}\,F_e\,,   \qquad q(0) = Q_{start}
\end{align}
where $q_{end} = q(T_e)$. The optimisation problem is now to find an optimal trajectory $q(t)$, satisfying operational constraints, that would maximise the target function~(\ref{eq:free_term_cond}). It is easy to see that any extremal trajectory of~(\ref{eq:free_term_cond}) is also extremal of~(\ref{eq:unperturbed_action1}, \ref{eq:unperturbed_action2}). Indeed, if the trajectory~$\bar q(t)$ is an extremal of~(\ref{eq:free_term_cond}), the variation of the target functional vanishes on any allowed small variation $\delta \bar q(t)$. Since there is no terminal boundary condition, the variation $\delta \bar q(t)$ may not vanish at the end point. But if the variation of the functional vanishes on the whole class of allowed variations $\delta\bar q(t)$, it also vanishes on the sub-class of variations $\delta \hat q(t)$ preserving terminal level ($\delta \hat q(T_e) = 0$).

We conclude that, in order to find a maximising trajectory for the functional~(\ref{eq:free_term_cond}), we need first to find a solution for the functional~(\ref{eq:unperturbed_action1}, \ref{eq:unperturbed_action2}) with fixed terminal condition. Then, considering this solution as a function of terminal state~$q_{end}$, we need to find a maximum of~(\ref{eq:free_term_cond}) as a plain function of~$q_{end}$, provided $q_{end}$ is not lying on the boundary. Thus, if the trajectory terminal point does not lie on the boundary, the target functional must satisfy the following condition
\begin{align}
    \frac{\p \tilde S}{\p q_{end}} = 0
\end{align}
This is valid for any target functional, regardless of whether the residual price is applied or not.

As shown in Sec.~\ref{sec:3}, if the target functional~(\ref{eq:unperturbed_action1}) is considered as a function of the terminal state~$q_{end}$, then its derivative with respect to $q_{end}$, in any interior point, is given by
\[
    \frac{\p S}{\p q_{end}} = -C
\]
Thus, the target functional~(\ref{eq:free_term_cond}), considered as a function of $q_{end}$, has a derivative
\begin{align}
    \label{eq:4.21}
    \frac{\p \tilde S}{\p q_{end}} = F_e - C
\end{align}
where $C$ is the trigger price of the final part of trajectory. We conclude that the condition for a target functional, with a free terminal state, to have a maximum is
\begin{align}
    \label{eq:4.22}
    C = F_e
\end{align}
which is valid if the terminal point is not lying on the boundary.

\subsubsection{terminal state cases}

We conclude that the terminal point of the optimal trajectory can only lie between the boundaries if $C=0$ (assuming no residual unit price). If $C>0$ then the trajectory inevitably terminates on the lower boundary. Similarly, if $C<0$, the terminal point lies on the upper boundary. The following relations must hold
\begin{align}
    \label{eq:term_level1}
    q_{end} = Q_{min}  \quad \Rightarrow & \quad C \geq 0 \\
    Q_{min} < q_{end} < Q_{max}  \quad \Rightarrow & \quad  C = 0 \\
    q_{end} = Q_{max} \quad \Rightarrow & \quad C \leq 0
\end{align}
The inverse is also true with small modifications
\begin{align}
    C > 0  \quad \Rightarrow & \quad q_{end} = Q_{min} \\
    C = 0  \quad \Rightarrow & \quad  Q_{min} \leq q_{end} \leq Q_{max} \\
    \label{eq:term_level2}
    C < 0  \quad \Rightarrow & \quad q_{end} = Q_{max}
\end{align}
If the residual unit price is paid, in the expressions~(\ref{eq:term_level1}-\ref{eq:term_level2}) we should substitute $C\to C-F_e$.

If the storage option has no residual unit price, it is most likely to terminate at the lowest possible volume level.  Indeed, the upper boundary of the dead zone cannot be lower than the smallest forward price:
\[
    C + \gamma_{rel} \geq \min(F(t))
\]
If the gas prices are positive, and operating costs are sufficiently small, then the trigger price is positive too. Hence the derivative of the option value with respect to the terminal level is negative, and thus the optimal trajectory must terminate at the lowest possible level. Consequently the storage option can be considered as a problem with fixed terminal state (although it may not be fixed according to the contract).

An opposite example is a \emph{swing} option. A typical swing option is a contract between a gas supplier and a trading company. The trader buys gas at some predefined strike price from the supplier and sells it on the market at the current market price. The spread between the market and the strike price becomes the effective gas price for the trader. The swing option allows the trader to take the gas from the supplier according to some flexible scheme (i.e. to decide, when to take and how much within certain constraints). Thus, the swing option can be formulated in terms of a storage option. However unlike the storage option, the effective price (the spread) in the case of a swing option can be both positive and negative. Depending on the actual prices and on the contract constraints, it may happen that the terminal state does not lie on either the lower or upper boundaries. In this case we deal with the problem with the free terminal condition. For this problem the trigger level has to be equal zero.

\subsection{Cycle constraint}
\label{sec:cycle_constraint}

A cycle constraint is sometimes imposed by the storage option seller in order to artificially restrict the option flexibility and reduce the option value. The cycle constraint is defined in~(\ref{eq:cycle_constrain}). For demonstration purposes we consider the intake cycle constraint. It caps the maximum injected volume during the operation period $t\in[0,T_e]$. Instead of considering it in the form of inequality we can proceed as follows. The solution is found in two steps. On the first step we solve the problem without the cycle constraint and calculate the cycle variable on the optimal trajectory. If the cycle variable is below the threshold $c_{max}$, the obtained solution satisfies the cycle constraint. However if the solution violates the cycle constraint, we can impose a ``biting'' cycle constraint\footnote{We say that the constraint is \emph{biting}, when an inequality ``$\leq$'' becomes an equality.}
\begin{align}
    \label{eq:cycle_constrain_eq}
    c(T_e) = \int_0^{T_e} \dot q(t)\, \theta( \dot q(t) )\,dt = c_{max}
\end{align}
This condition allows us to formulate our problem in terms of conditional variational extrema. It is solved by means of a Lagrange multiplier, namely, we introduce a modified Lagrangian
\begin{align}
    \label{eq:modif_lagr_cycle}
    L = - \Big[ \dot q F + \phi(q) + \psi(\dot q) + \gamma(\dot q) + \lambda\, \dot q(t)\, \theta( \dot q(t) ) \Big]
\end{align}
where $\lambda$ is an independent variable (the Lagrange multiplier). The solution to the latter problem will contain the undefined coefficient $\lambda$, which can be found from the additional equation~(\ref{eq:cycle_constrain_eq}).

Within the boundaries, the Lagrange-Euler equation for the constrained Lagrangian becomes (compare with Eq.~(\ref{eq:within_boundaries}))
\begin{align}
    C & = F(t) + \psi'(\dot q) + \gamma'(\dot q)  + 
        \lambda \,\Big(\theta(\dot q) + \dot q\,\delta(\dot q) \Big)  \nonumber \\
    & = F(t) + \psi'(\dot q) + \gamma'(\dot q)  + \lambda \,\theta(\dot q)
\end{align}
Solving this equation we find that the solution is similar to the solution without the cycle constraint, but has a ``dead zone'' which is wider by an amount $\lambda$. Thus, the cycle constraint is equivalent to having additional operating costs. One can find fictive additional release and/or injection costs (constant for the whole trajectory), such that, if added to the original storage, will result in an optimal trajectory respecting the cycle constraint. It is inconsequential if the additional cost is for injection, release, or both, since it is only the width of the dead zone which matters. The only exception to this rule is if the optimal trajectory does not inject or release enough volume. For instance, if the trajectory never releases gas then the additional release cost will not effect the trajectory, and thus will not help respect the intake cycle constraint. Practically, one could use additional injection costs to fulfil the intake cycle constraint, and release costs to fulfil the release cycle constraint.

Generally the solution between boundary touches contains two undetermined parameters -- the trigger price $C$ and the cycle Lagrange multiplier $\lambda$.  They can equivalently be expressed by the upper and lower bounds of the dead zone. These parameters have to be found in such a way that the obtained solution satisfies the boundary touch conditions~(\ref{eq:bc1}-\ref{eq:bc4}), initial and terminal conditions~(\ref{eq:init_term_cond}), and the cycle constraint~(\ref{eq:cycle_constrain_eq}).

Notice that the additional operating cost is the same for the entire time horizon, regardless of whether the trajectory touches the boundary or not.

This result provides an efficient way of numerical computation of a storage problem with a cycle constraint. First one should solve the problem without any modification, and check the cycle variable. If the original solution breaches the cycle constraint, one needs to find a proper virtual operating costs, which allow the cycle constraint to be met. The search for the right virtual operating cost can be done by a ``shooting'' method, and in principle, should converge after a few iterations, since the cycle variable is obviously is monotonically decreasing function of the operating costs.

\subsection{Volume dependent injection/release rates}
\label{sec:volume_dependent_rates}

If the maximal injection/release rates are volume dependent, the penalty function $\psi$ becomes an explicit function of volume: $\psi = \psi(q,\dot q)$. In this case the Euler-Lagrange equation within the boundaries becomes
\begin{align}
    \label{eq:4.25}
    \p_q \,\psi(q,\dot q)  = \frac{d}{dt}\Big[F + \p_{\dot q}\,\psi(q,\dot q) + \gamma'(\dot q) \Big]
\end{align}
The l.h.s. of this equation is an integrable implicit function of time. Designating
\[
    C(t) = F + \p_{\dot q}\,\psi(q,\dot q) + \gamma'(\dot q)
\]
we get
\[ 
    \frac{dC(t)}{dt} = \p_q \,\psi(q(\tau),\dot q(\tau)) \,;\qquad C(t) = C_0 + \int_0^t  \p_q \,\psi(q(\tau),\dot q(\tau)) \,d\tau 
\]
We thus obtain the same solution as in the case of volume independent rates with the difference that $C$ is not constant but a function of time. We conclude that the optimal exercise strategy still preserves the ``bang-bang'' property. However, it does not possess a constant trigger level. The function $C(t)$ can be reverse-engineered from the optimal trajectory $q(t)$, provided the latter is known. 
%We show how to calculate $C(t)$ in appendix~\ref{sec:trigger_level_vol_dependent}.

\section{Optimal stochastic exercise}
\label{sec:StochTriggerPrice}

In this paper we do not consider the stochastic problem in detail. However, we do wish to address one particular aspect of the storage option problem in its stochastic formulation, namely, the optimal stochastic exercise.

In the stochastic formulation the prices are considered to be randomly changing during trading, as opposed to being fixed in the intrinsic problem formulation. The forward curve now becomes a function of two parameters -- observation time $t$ and delivery time $T$:  $F(t,T)$. For every fixed delivery time $T$ the forward curve becomes a stochastic process indexed by the observation time $t$. A choice of the price process poses a problem per se, and to some extent is driven by knowledge of the analytic solution properties.

One peculiarity about the stochastic problem is that the optimal exercise profile can not be found upfront for the entire time horizon. It has to be continuously adjusted to the changing market prices. However for the current observation time $t$ there will always exist an optimal stochastic exercise decision $\dot q_{st}(t)$, which remains valid for an infinitesimal time interval $dt$. The optimal exercise is conditional on the current state of the market $F(t,T)$ and of the storage $q(t)$, or expressed in probabilistic language, conditional on the filtration $\mathcal{F}_t$. We refer to the optimal exercise decision at time $t$ as the optimal \emph{prompt exercise} (the prompt exercise is a continuous time analogue of the spot trade in the discrete time systems).

In a similar way we can introduce an \emph{intrinsic prompt exercise}. For every time $t$ we can formulate a new intrinsic problem, for which $q(t)$ plays a role of the initial condition, and $F(t,T)$ is the forward curve (as a function of $T$). The solution of this problem is an optimal intrinsic exercise $\dot q_{in}(t,T)$, for which observation time $t$ plays a role of a parameter, and the time derivative is taken with respect to $T$. The intrinsic prompt exercise is defined as the initial value of the optimal exercise profile $\dot q_{in}(t,t)$. Below we use the notation 
\[
    r(t) := \dot q_{in}(t,t)
\]
for the intrinsic prompt exercise.

In this section we establish a link between the optimal stochastic $\dot q_{st}(t)$ and intrinsic prompt exercise $r(t)$, and show why the intrinsic exercise is a good approximation for the stochastic one. We do not consider injection/release costs for the sake of simplicity of notation. The generalisation is straightforward.

Although the main result of the current section turns out to be independent of the choice of price process, for the sake of completeness we discuss some of most common approaches to modelling the stochastic price evolution.

\subsection{Forward curve evolution model}

Usually the forward price evolution is approximated by a zero mean stochastic process governed by the stochastic differential equation of the form
\begin{align}
    \label{eq:generic_price_process}
    &\frac{dF(t,T)}{F(t,T)} = dM_t(T)
\end{align}
where $M_t(T)$ is a martingale process of time $t$, parametrised by the variable~$T$. For the derivation of the optimal stochastic exercise the particular choice of the price process is not essential. The only important assumption about the price process is that it is drift-less, i.e. the mean of the price increment (for every delivery time $T$) is vanishing:
\begin{align}
    \la dF(t,T) \ra = 0
\end{align}
This assumption is plausible on liquid markets without friction.

For different delivery times $T$ the processes $M_t(T)$ may or may not be correlated. One particular example of the price process is the 1-factor model
\begin{align}
    \label{eq:1-factor_proc}
    \frac{dF(t,T)}{F(t,T)} = \sigma\,e^{-\alpha (T-t)}\,dW_t
\end{align}
where $W$ is the standard Wiener process. The model was first presented in \cite{schwartz1998valuing}, its application to the Energy market can be found in \cite{clewlow1999valuing}. The one factor model is very popular due to its simplicity. Notice that the processes $M(T)$ for all delivery times $T$ are perfectly correlated.

\subsection{Stochastic trigger price}
\label{sec:4.6.4}

Let $F(t,T)$ be the forward curve observed at time $t$. We define the \emph{spot price} process $s(t)$ as
\begin{align}
    s(t) := F(t,t)
\end{align}
Since the forward price curve evolves with time, the spot price $s(t) = F(t,t)$ generally does not follow the original forward curve ($s(t) \neq F(0,t)$ almost everywhere). The spot process $s(t)$ is uniquely defined by the forward price process.

Consider a set of possible spot price curves $s_i(t)$, each having corresponding probability $p_i$. Viewing a particular spot price curve $s_i(t)$ as a deterministic function, we can solve a static optimisation problem for the target functional\footnote{The discussion of why the expressions for the target functional are well-defined would lead us too far from the scope of the paper. However we should point out that the problem can be formulated in a discrete time, where it becomes a well-defined finite-dimensional problem. The transition to the continuous time can be done if we demand that the spot price process is sufficiently ``good''. Wiener process is known to be almost surely continuous, and the integral of Wiener process is well-defined.}
\[
    S_i = -\int_0^{T_e} \dot q(\tau)\,s_i(\tau)\,d\tau
\]
defined for the time moment $t=0$. For every spot price curve $s_i(t)$ we can find a corresponding trigger price $C_i$, optimal exercise strategy $\dot{ \bar q}_i(t)$ and intrinsic value $S_i$. The prompt intrinsic exercise
\[
    r_i := \dot{ \bar q}_i(0)
\]
for the time moment $t=0$ is optimal on the $i$-th spot price path.

Now we are searching for an optimal prompt exercise decision $\dot q_{st}(0)$ for the stochastic problem at time $t=0$. We omit the argument $0$ in what follows to simplify the notation. The set of the spot price paths with corresponding probabilities contains the entirety of the stochastic information required for making an optimal stochastic exercise decision. Since we can only make one exercise decision at a time, this exercise decision will not be optimal for some spot price paths. If on the $i$th spot price path the optimal exercise is $r_i$, the sub-optimal exercise decision $\dot q_{st}\neq r_i$ will lead to a loss of value on that path. We need to make such a choice of $\dot q_{st}$ that would minimise the expected losses over all paths for which the decision is not optimal.

During the time interval $dt$ the volume increment in the storage is given by $dq = \dot q_{st}\,dt$. Thus, on the $i$th spot price path the volume in the storage is by $(\dot q_{st} - r_i)\,dt$ bigger than it would be if the exercise was optimal for that path. This leads to a change of the value on the $i$th path by
\begin{align}
    \label{eq:4.13}
    \delta S_i = \Big( C_i - F(0,0) \Big)\Big( \dot q_{st} - r_i\Big)dt
\end{align}
Indeed, the change of the value occurs for two reasons: the additional expense $\delta S_{i1}$ of buying additional volume, and an increase in value $\delta S_{i2}$ due to a change in the actual storage volume. The price of additional volume is
\[
    \delta S_{i1} = - F(0,0)\, \Big( \dot q_{st} - r_i\Big)\,dt
\]
The change in value due to a change in the initial volume follows from the formula $\p S/ \p q_{start} = C$ (see Sec.~\ref{sec:3}):
\[
    \delta S_{i2} = C_i \,\Big( \dot q_{st} - r_i\Big)\,dt
\]
Combining the latter two equations we obtain Eq.~(\ref{eq:4.13}). Now, averaging $\delta S$ over all paths, we have
\begin{align}
    \la \delta S \ra = \sum_i p_i \Big( C_i - F(0,0) \Big)
        \Big( \dot q_{st} - r_i\Big)dt
\end{align}
where $p_i$ is the probability of the $i$th path. The optimal exercise decision $\dot q_{st}$ should maximise the expected change of option value $\la \delta S \ra$. Deriving $\la \delta S \ra$ with respect to $\dot q_{st}$ we obtain
\begin{align}
    \frac{\p \la \delta S \ra}{\p \dot q_{st}} = \sum_i p_i\,\Big( C_i - F(0,0) \Big)\,dt =
        \Big( \la C \ra  - F(0,0)\Big) \,dt
\end{align}
where $\la C \ra = \sum_i\,p_i\,C_i$ is the trigger price averaged over all possible spot price processes. If the spot price $F(0,0)$ is bigger than the average trigger price $\la C\ra$, then the derivative is negative, and the optimal stochastic exercise is the smallest possible allowed by the constraints. Similarly if $F(0,0)$ is smaller than the average trigger price, then the optimal exercise is the biggest possible allowed by the constraints. Thus, the optimal prompt stochastic exercise is
\begin{align}
    \dot q_{st} = \left\{
        \begin{array}{ll}
            r_{min} \qquad &  F(0,0) > \la C \ra              \\
            r_{max} \qquad &  F(0,0) < \la C \ra            
        \end{array}
    \right.
\end{align}
We see that the optimal stochastic exercise is again bang-bang, and that the expected value of the intrinsic trigger price $\la C \ra$ can be interpreted as the stochastic trigger price $C_{st}$:
\begin{align}
    \label{eq:stochastic_trigger_price_1}
    C_{st} = \la C \ra
\end{align}

Next we show that, under some assumptions, in the leading order of the Taylor expansion the stochastic trigger price equals the intrinsic trigger price. Indeed, let us designate $s_0(t)$ as the spot price path which coincides with the initial forward curve at time $t= 0$:
\[
    s_0(t) \equiv F(0,t)
\]
The trigger price $C_0$ is then the intrinsic trigger price, calculated at time $t=0$ on the basis of the forward curve $F(0,t)$. All other paths $s_i(t)$ are different, and we designate
\[
    \delta s_i(t) = s_i(t) - s_0(t)
\]
This difference is zero at the beginning of time period: $\delta s_i(0) = 0$ since all paths of the spot price process start in the same point $s_i(0) = F(0,0)$. At any future point in time the difference $\delta s_i(t)$ may be arbitrary large, but the majority of paths stay within the standard deviation
\[\sqrt{\lla \delta s^2(t) \rra}\]
If $\delta s_i(t)$ remains small, then the difference in the trigger level $\delta C_i = C_i-C_0$ can be calculated by means of the first order Taylor expansion~(see App.~\ref{sec:variation_trigger_price}):
\[
    \delta C_i \approx \frac{\Delta r}{K} \int_0^{T_e} \delta[C_0-s_0(t)]\,\delta s_i(t)\,dt
\]
where we designated
\begin{align*}
    & \Delta r = r_{max} - r_{min} \qquad 
    K = \Delta r \,\int_0^{T_e} \delta[C_0 - s_0(t)]\,dt
\end{align*}

Since the forward price process is a martingale, we conclude that $\la \delta s_i(t) \ra = 0$, and hence in the first order $\la \delta C_i \ra = 0$. Consequently
\begin{align}
    \label{eq:stochastic_trigger_price}
    C_{st} = \la C \ra = C_0 + \la \delta C \ra \approx C_0
\end{align}
which means that the stochastic trigger price equals the intrinsic trigger price in the first order of the Taylor expansion.

We performed the calculation of the stochastic trigger price for the initial time moment $t=0$. For any later time we can repeat the same calculation, and thus the following statement holds: for any time $t$ the optimal stochastic exercise is bang-bang, and the stochastic trigger price $C_{st}(t)$ is approximately equal to the intrinsic trigger price $C(t)$, calculated at time $t$ on the basis of the deterministic forward curve $F(t,T)$ and initial condition $q(t)$. The exact stochastic trigger price can be found as an average over the intrinsic trigger prices calculated for the entire set of spot price paths having $F(t,T)$ as initial condition.

Note that to derive the relation (\ref{eq:stochastic_trigger_price_1}) we did not use any approximation, and hence this relation is exact. In order to derive relation~(\ref{eq:stochastic_trigger_price}) we used only the assumption of small variation
\begin{align}
    \label{eq:condition_for_rolling intrinsic}
    \frac{ \sqrt{\lla \delta s^2(t) \rra}}{s_0(t)}\ll 1
\end{align}
For the standard 1-factor price process~(\ref{eq:1-factor_proc}) the relative spot price variation can be easily estimated,
\[
    \frac{ \sqrt{\lla \delta s^2(t) \rra}}{s_0(t)} =
    \sqrt{ \exp \left[ \frac{\sigma_0^2}{2\,\alpha}\,\Big(1- e^{-2\,\alpha\,t}\Big)\right]-1 }
\]
In particular, for $\alpha\,T_e\ll 1$ the condition for the intrinsic exercise approximation becomes
\[\sigma_0\,\sqrt{ T_e} \ll 1
\]
and for $\alpha\,T_e\gg 1$ we obtain the condition
\[ \frac{\sigma_0}{\sqrt{2\,\alpha}} \ll 1 
\]

\subsubsection{Equivalent stochastic trigger prices}

The intrinsic trigger price  is a well-defined value, which has an impact on the whole intrinsic exercise strategy (or at least for the whole part of trajectory between boundary touches). One can interpret the trigger price as an alternative way of representing the exercise strategy, since the optimal exercise can be derived for any point in time from the trigger level and price.

To the contrary, the stochastic trigger price only has an impact on the prompt exercise decision over the infinitesimal time interval $dt$. The trigger price may be different for every observation time, since, due to the stochastic nature of the prices, there is no optimal exercise strategy that can be calculated for the whole time horizon of the storage contract.

The stochastic trigger price is thus not uniquely defined. Indeed, if some price $C(t)$ is a trigger price at time $t$, then it defines the prompt exercise according to the rule: inject if the spot price $s(t)$ is below $C(t)$, and release if the spot price is above $C(t)$. In the first case, any price $\tilde C(t)$ which is above the spot price would lead to the same exercise decision, and thus play a role of the trigger price. Similarly, in the second case, any price $\tilde C(t)<s(t)$ would lead to the same prompt exercise. Thus we can say that for each time moment there exists a class of \emph{equivalent} stochastic trigger prices.

If $C_1(t)$ and $C_2(t)$ are two equivalent stochastic trigger prices, and $s(t)$ is the spot price, then obviously
\begin{align}
    & s(t)<C_1(t) \quad \Leftrightarrow \quad s(t)<C_2(t)\\
    & s(t)>C_1(t) \quad \Leftrightarrow \quad s(t)>C_2(t)
\end{align}
Form these statements we conclude that there is one case when all equivalent trigger prices should coincide
\begin{align}
    \label{eq:coinciding_trigger_prices}
    &  s(t)= C_1(t) \quad \Leftrightarrow \quad s(t)=C_2(t)
\end{align}

\subsubsection{Value loss due to sub-optimal exercise}

Next we show that, similar to intrinsic case, the stochastic trigger price (one from the range of equivalent prices) is related to the stochastic option value $V$ as
\begin{align}
    \label{eq:stoch_trigger_price_derivative}
    C_{st} = \frac{\p V}{\p q}
\end{align}
where $q$ is the volume level. 

Let $V(q,t)$ be the stochastic option value at time $t$, considered as a function of the volume in storage $q$ and time $t$ (for a given forward curve $F(t,T)$). Let the spot price be $s = F(t,t)$. We consider a portfolio  composed of the storage option (over the remaining time horizon), hedge trades and cumulative cash-flow generated by the spot exercise trades. The increment of the portfolio value $d\Pi$ over the time interval $dt$ is given by
\begin{align}
    d\Pi = -s\,dq + \frac{\p V}{\p q}\,dq + \ldots % \frac{\p V}{\p t}\,dt
\end{align}
The first term is due to the cash flow from purchasing volume $dq$ at the spot price $s$, and the second is due to changed option value as a function of volume $q$. We omitted the terms originating from derivatives with respect to time and forward price, as they have no impact on the derivation of the trigger price. The portfolio value increment is maximised by the following choice of volume increment
\begin{align}
    dq = \left\{
        \begin{array}{ll}
            dq_{max} = r_{max}\,dt \qquad &  s < \p V /\p q  \\
            dq_{min} = r_{min}\,dt \qquad &  s > \p V /\p q 
        \end{array}
    \right.
\end{align}
We see that the derivative $\p V/\p q$ can be interpreted as a stochastic trigger price. It is equivalent to the trigger price $C_{st}$ calculated previously. 

The trigger price allows us to estimate a loss of portfolio value due to a sub-optimal exercise. Considering the portfolio increment $d\Pi$ as a function of the exercise $dq$ we can write
\begin{align}
    \label{eq:portfolio_variation}
    \frac{\p \,d\Pi}{\p\, dq} = C_{st} - s 
\end{align}
Let $dq$ be an optimal exercise, and $d\tilde q$ be some arbitrary (sub-optimal) exercise, which leads to the portfolio increment $d\tilde \Pi$. Designating $\delta q = d\tilde q - dq$ as the deviation of the prompt exercise from the optimal one, and $\delta \Pi = d\tilde\Pi - d\Pi$ as the corresponding value loss, we obtain
\begin{align}
    \label{eq:portfolio_loss}
    \delta \Pi = \left(C_{st}(t) - s(t)\right)\,\delta q
\end{align}
Notice that the portfolio value is insensitive to the sub-optimal exercise at the trigger times $t^*$ defined by the equation $C_{st}(t^*) = s(t^*)$.

\section{Discussion and numerical analysis}
\label{sec:discussion}

In this paper we presented a solution to the intrinsic storage optimisation problem. It is represented in an implicit form, as the solution contains an undefined parameter -- the trigger price -- which can not be expressed analytically. However, the solution gives us a deep insight into the essence of the storage problem, and provides evidence for sometimes very counter-intuitive features of the storage option.

Here we first give an overview of the formal solution, and then discuss the most important conclusions that can be drawn from the analytic results.

\begin{itemize}
    \item An optimal exercise of the intrinsic problem can be expressed by means of the trigger price. This is such a price level $C$ for which the optimal exercise decision can be represented as
    \begin{align}
        \dot q(t) = \left\{ \begin{array}{ll}
                        r_{max}(t)\qquad & F(t) < C \\
                        r_{min}(t)\qquad & F(t) > C
                    \end{array}
                    \right.
    \end{align}
    If injection and release are subject to operating costs then there is also a ``dead zone'' around the trigger price. The optimal exercise in this case is represented in Eq.~(\ref{eq:solution}). It contains an additional ``do nothing'' state, which is triggered when the price is within the dead zone. The width of the dead zone is solely defined by the operating costs.
    
    In the more general case, an optimal exercise trajectory may partially go along the upper or lower boundary. The entire trajectory then consists of several parts of ``active'' injection/release, separated by a boundary touch. The trigger prices for different ``active'' parts separated by a boundary touch may be different.
    \item If the optimal exercise trajectory $q(t)$ touches the storage boundary (i.e. reaches the storage minimum or maximum volume), then conditions~(\ref{eq:bc1}-\ref{eq:bc4}) must be satisfied. In the case of no operating costs these conditions are simplified to
    \begin{align}
        F(t^*) = C
    \end{align}
    where $t^*$ is the time of the boundary touch.
    \item If the option value $V$ is considered as a function of the terminal level $q_{end}$, then is has a simple relation to the trigger price given by
    \begin{align}
        \frac{\p V}{\p q_{end}} = F_e - C
    \end{align}
    where $F_e$ is a residual unit price (price paid back for the gas remaining in the storage after the end of the contract).
\end{itemize}

Now we are going to discuss some of the most significant practical implications, which follow from our analysis.

\begin{itemize}    

    \item An important statements can be made about the terminal level of the storage. First consider a storage option with no residual unit price. If the prices are strictly positive and if operating costs are sufficiently small, then the trigger price is positive, and hence the derivative $\p V/\p q_{end}$ is negative. Thus the optimal trajectory has to terminate at the lowest possible storage level
    \begin{align}
        q_{end} = Q_{min}
    \end{align}
    With a good accuracy the same applies to the stochastic storage option.
   %\footnote{Strictly speaking, in a stochastic formulation there is always a chance that the biggest difference between the forward prices becomes smaller than the width of the dead zone. In this case the trigger price can become non-positive, and hence trajectory may terminate above $Q_{min}$. However in realistic situation the probability of such event is extremely low.}. 
    
    \item Unlike storage, for a swing option it is possible to terminate between global volume constraints. Indeed, an effective price of the underlying of a swing contract is a spread between the market price and the strike price, which can be positive or negative. Hence the sign of the trigger price is indefinite. The actual terminal volume (in other words, the total taken volume) depends on the contract ``moneyness''. If the swing contract is deep in the money (spread is positive), the trigger price is positive and the exercise trajectory will terminate on the lower boundary. If the contract is deep ``out of the money'', the spread is negative, and so is the trigger price. Hence the exercise trajectory will end on the upper boundary. For an ``at the money'' option the total taken volume might lie between the global constraints. In this case the trigger price must vanish
    \begin{align}
        C = 0
    \end{align}
    The same applies to the storage option, if a terminal unit price is paid for the remaining volume. Depending on the moneyness, the terminal volume may lie between the volume constraints.
    
    Notice that a vanishing trigger price implies a simple exercise strategy: release as soon as market price is above the strike. In this case the swing contract becomes a strip of independent vanilla call options. If the terminal point of the intrinsic exercise trajectory is sufficiently far from the global constraints, then the same applies to the stochastic swing option: it can be priced as a strip of \emph{independent} vanilla call options. The independence follows from the fact that the majority of stochastic exercise trajectories will also terminate between the boundaries. Hence the trigger price $C = 0$ is preserved, and the exercise rule for all maturities remains the same. A swing option which can be represented as a strip of independent options has a significant advantage, since it can be replicated (and hence priced and hedged) by the market observables (of course, only to the extent allowed by the market liquidity).
    
    \item If a carry cost is imposed, then the trigger price becomes a growing function of time satisfying
    \begin{align}
        \frac{ dC(t)}{dt} = \gamma_c(t)
    \end{align}
    where $\gamma_c(t)$ is the carry cost per unit volume per unit time. As has been shown in the main text, in terms of exercise strategy, carry cost is equivalent to an additional discount factor
    \[
        D(t) = 1 - \frac{1}{F(t)}\int_0^t \gamma_c(u)\,du \approx
        \exp\left(- \frac{1}{F(t)}\int_0^t \gamma_c(u)\,du\right)
    \]
    Notice that the discount factor leads to a correct exercise strategy, but wrong cash flow. The correct cash flow has to be calculated after the optimal exercise is found.
    
    \item A cycle constraint is equivalent to an effective additional injection/release cost. Practically this means that in order to meet the cycle constraint one can use the following iterative approach. First solve the intrinsic problem without any modification. If the cycle constraint is breached, one should add an operating cost (let us call it virtual opex) and solve the new modified problem. 
    Too low a virtual opex might not be sufficient to meet the cycle constraint, too high a virtual opex will reduce the cycle variable below the constraint. After a few iterations one should converge to an optimal choice of virtual opex, which leads to a biting cycle constraint. Such an iterative solution is computationally much more efficient than a brute force approach, since the latter would require an additional dimension in the state space.
    
    Notice that, similar to the case with carry cost, the virtual opex gives rise to a correct exercise strategy, but wrong cash flow. The correct cash flow has to be calculated after the optimal exercise is found.
    
    \item One of the most important results of our analysis concerns the stochastic spot exercise. We have shown that the optimal stochastic exercise is ``bang-bang'' in exactly the same way as it is for the intrinsic problem. Thus the stochastic exercise can be expressed in terms of the \emph{stochastic trigger price} $C_{st}(t)$. Unlike the intrinsic problem though, the stochastic trigger price cannot be defined upfront for the entire time horizon. 
    For each time moment $t$ it can only be determined for a short time interval, over which the forward curve does not change significantly as a function of the observation time. Typically gas storage and swing contracts are exercised with a daily granularity (day-ahead nominations), and hence one day is a natural time interval of the trigger price validity. A more interesting use case for the trigger price is posed by weekends and bank holidays, when there is no gas market activity, and the trigger price should be prolonged over a period of a few days.

    \item The stochastic trigger price turns out to be very close to the intrinsic trigger price, or equivalently, the stochastic prompt exercise is almost identical to the intrinsic prompt exercise. The importance of this fact is difficult to overstate, as it justifies the so called rolling intrinsic (RI) strategy. This strategy is a rule to exercise the storage option using the intrinsic optimisation on every time step. 
    The RI strategy is performed as follows. At each time step $t$, given the current market prices $F(t,T)$ and the current storage level $q(t)$, one solves the intrinsic problem and finds the optimal intrinsic exercise profile $\dot q(t,T)$ for the remaining time horizon. The prompt exercise (the exercise trade) $\dot q(t,t)$ is traded on the spot market, and the remainder (the hedge trade) on the forward market (of course, when performing exercise and hedge trades, one has to take into account the existing hedge position). This action is repeated each trading day until the end of the storage contract.
    
    The rolling intrinsic strategy is one of the most frequently used by traders. This strategy is attractive for several reasons, one of which is that it is trivial to follow. In particular, it does not require any knowledge of, or calibration to, the price process. 

    Since the intrinsic hedge is not perfect (as it does not take the stochastic price dynamics into account), the terminal portfolio value is to some extent random, and may be significantly different for different price evolution scenarios even with identical initial conditions.
    An important fact is that, on average, the terminal value of the portfolio is only determined by the prompt exercise\footnote{Indeed, the hedge trades consist of pure financial instruments. Since the price process is a martingale, the average value of all hedge contracts vanishes. It is only the ``physical'' prompt exercise trades which generate a non-trivial average cash flow.}. 
    Because the intrinsic prompt exercise, as we have seen, is very close to the correct stochastic prompt exercise, we conclude that the RI strategy on average generates the correct terminal option value (in a risk neutral probability measure).

    \item A correct stochastic option value can only be achieved by following an optimal stochastic exercise strategy. However there are a few time moments during the storage contract life time, when the prompt exercise decision is irrelevant. Indeed, according to Eq.~(\ref{eq:portfolio_loss}) the deviation of a portfolio value from its optimum is given by
    \begin{align}
     \delta \Pi = (C_{st}(t) - s(t))\,\delta q
    \end{align}
    where $\delta q$ is the deviation from the optimal exercise, and $s(t)$ is the spot price. We see that the portfolio value loss $\delta \Pi$ becomes insensitive to the wrong exercise every time the spot price crosses the trigger level, i.e. at trigger times.
    
    An important consequence follows from this statement. Since the stochastic and intrinsic trigger prices are very close, they may only lead to a different exercise decision at those time moments when the spot price gets close to the trigger price, in other words -- at trigger times. 
    However, if instead of a correct stochastic exercise one performs an intrinsic exercise, it won't lead to a change of option value, since the option value is insensitive to the exercise decision at the trigger times. This gives yet another justification why the intrinsic exercise is a good approximation of a correct stochastic exercise.

    Another interesting consequence is that, according to the boundary touch conditions~(\ref{eq:bc1}-\ref{eq:bc4}), every time the trajectory is approaching the boundary, the option value becomes insensitive to the exercise decision.
    
\end{itemize}

We performed a numerical simulation in order to demonstrate the effect of the intrinsic exercise on the PnL. We simulated a PnL evolution of a portfolio consisting of a 1-year storage option, a linear hedge and a cumulative cash flow. For the simulation we generated 300 different evolution scenarios of a forward curve, driven by the one-factor price process~(\ref{eq:1-factor_proc}) with parameters $\alpha=20$ and $\sigma = 0.9$. 

For each price evolution scenario we ran two instances of the portfolio, one of them using a correct stochastic exercise decision as a benchmark, another using an exercise decision based on intrinsic valuation at every time step. 

Since the storage pricing model used the same assumptions as those used for simulating a forward curve, the first portfolio is perfectly hedged, and hence its PnL follows a straight line plus some noise which is due to the discrete time and Monte Carlo valuation error. 

The second portfolio used correct stochastic hedging, but intrinsic prompt exercise. The PnL of the second portfolio may be slightly different only due to the rare events when intrinsic and stochastic exercise decisions are different. Our aim was to demonstrate that the PnL of the two portfolios is essentially the same.

\begin{table}[htb]
    \centering
        \begin{tabular}{l|c}
            \hline
            \hline
            Full option value & $305.8\pm 0.5\$$ \\ \hline
            Time value & $31.1 \pm 0.5\$$ \\ \hline
            Terminal PnL (stochastic exercise) & $306 \pm 0.2\$$ \\ \hline
            Terminal PnL (intrinsic exercise) & $305.4\pm 0.2\$$ \\ \hline
            PnL loss (intrinsic exercise) & $0.6 \pm 0.1\$$ (2\% of the time value) \\ \hline
        \end{tabular}
    \caption{\small Simulation results: intrinsic trigger price yields correct exercise decision.}
    \label{tab:pnl}
\end{table}

The results of the simulation are presented in the table~\ref{tab:pnl} (we used the standard error for the accuracy interval). We found that the intrinsic exercise led to a systematic PnL loss of about 2\% of the option time value. 

For our example we deliberately chose such parameters (fast storage and fast oscillating forward curve) to enforce a large number of trigger times, in order to make the effect of intrinsic exercise more pronounced. More realistic examples of the storage option would have much fewer trigger times, and hence even smaller PnL loss. Thus the simulation confirmed that the intrinsic exercise has a negligible effect on the PnL.

\paragraph{acknowledgements}
I'm grateful to my colleague Garry Bowen for reviewing the manuscript and giving his valuable comments.

\vspace{0.5cm}

\begin{flushleft}
{\Large \bf Appendices}
\end{flushleft}

\appendix

\appendix\normalsize

% \section*{Appendices}
% \addcontentsline{toc}{section}{Appendices}
% \renewcommand{\thesubsection}{\Alph{subsection}}

\section{Modified and unmodified functionals}
\label{sec:modified_unmodified_functionals}

Let $\bar q(t)$ be the optimal solution of the modified functional $S[q(t)]$. Here we show why $\bar q(t)$ respects the constraints of the original problem, and why $\bar q(t)$ is also an optimal solution for the unmodified functional $S_0[\bar q(t)]$ with constraints.

The optimal trajectory $\bar q$ belongs to the broad class of locally integrable functions $\mathcal{L}_1$. Let us introduce a class of trajectories $\mathcal{K}\subset\mathcal{L}_1$ satisfying all the constraints. By definition, penalty functions vanish on $\mathcal{K}$:
\[
    \phi(q)\equiv 0 \qquad \psi(\dot q)\equiv 0 \qquad \text{iff}\quad q\in\mathcal{K}
\]
Obviously the values of modified and unmodified functionals coincide on $\mathcal{K}$:
\[
    S[q] = S_0[q] \qquad \forall q\in\mathcal{K}
\]

Next we show that the optimal trajectory $\bar q$ belongs to the class $\mathcal{K}$ and hence satisfies the constraints. Indeed, generally on some parts of the extremal trajectory $\psi'(\dot{\bar q}) \not = 0$  and $\phi'(\bar q) \not = 0$. According to Eqs.(\ref{eq:within_boundaries}) and (\ref{eq:on_the_boundary}) they take  finite values. One can show (we leave it without proof) that in the limit $N_{\phi,\psi}\to\infty$ everywhere where $\phi'$ (or $\psi'$) take finite values, the function $\phi$ (or $\psi$) vanishes. We conclude that on the extremal trajectory
\begin{align}
    \label{eq:zero}
    \phi(\bar q) \equiv 0 \qquad
    \psi(\dot{\bar q}) \equiv 0
\end{align}
and hence $\bar q \in\mathcal{K}$. Form~(\ref{eq:zero}) follows that the values of modified and unmodified target functionals coincide on the maximal trajectory:
\[S[\bar q] = S_0[\bar q]\]

The last step is to prove that the trajectory $\bar q$ delivers the maximum for the unmodified functional on the space~$\mathcal{K}$. It is obvious from the following consideration. The trajectory $\bar q$ maximises the modified functional on the space $\mathcal{L}_1$ and hence it will also maximise this functional on the smaller space~$\mathcal{K}\subset\mathcal{L}_1$. Since the values of both functionals coincide on~$\mathcal{K}$, the trajectory $\bar q$ will also maximise the unmodified functional $S_0$ on the space of all trajectories satisfying the constraints.

%
%
%\section{Trigger level for volume dependent rates.}
%\label{sec:trigger_level_vol_dependent}
%
%
%Let $q(t)$ be the optimal trajectory. According to Eq.~(\ref{eq:4.25}) on this trajectory the derivative $\p_{\dot q}\,\psi(q,\dot q)$ takes a finite value. This value is related to the trigger price $C(t)$ as
%\[
    %\p_{\dot q}\,\psi(q(t),\dot q(t)) = C(t) -F(t) - \gamma'(\dot q(t))
%\]
%We can also find the relation between the derivatives $\p_q\,\psi(q,\dot q)$ and $\p_{\dot q}\,\psi(q,\dot q)$. Indeed, let us consider a ``phase'' space $\{q,\dot q\}$. The optimal trajectory in this space is represented in a parametric form $\{q(t), \dot q(t)\}$. If the trajectory $q(t)$ lies, for instance, on the maximum injection rate boundary $\dot q = r_{max}(q)$, then the relation between $dq$ and $d\dot q$ is simply $d \dot q = r'_{max}(q)\,dq$. It's easy to show that the gradient $\{\p_q\,\psi(q,\dot q),\p_{\dot q}\,\psi(q,\dot q)\}$ must be orthogonal to the vector $\{dq, d\dot q\}$. Hence
%\[
    %\p_q\,\psi(q,\dot q) = r'_{max}(q)\,\p_{\dot q}\,\psi(q,\dot q)
%\]
%We thus obtain the differential equation
%\[
    %\frac{dC(t)}{dt} = r'(q(t)) \Big( C(t) -F(t) - \gamma'(\dot q(t)) \Big) \qquad\quad
    %r'(q(t)) = \left.\frac{dr(q)}{dq}\right|_{q = q(t)}
 %\]
%Here $r = r_{max}$ must be used for ``injection'' and $r=r_{min}$ for ``release'' parts of the trajectory. This equation allows us to find the function $C(t)$ provided the optimal trajectory $q(t)$ is already known.
%

\section{Variation of the trigger price}
\label{sec:variation_trigger_price}

Here we calculate how the intrinsic trigger price is modified by a small variation of the forward curve. We only estimate the first order Taylor expansion in this section. 

Let $\delta F(t)$ be the variation of the forward curve, and $\delta C$ is the corresponding variation of the trigger price. In order to derive the relation between $\delta F(t)$ and $\delta C$, we depart from the relation connecting the initial and terminal conditions:
\begin{align}
    \label{eq:boundary_condition2}
    Q_{end} = Q_{start} + \int_0^{T_e} \dot {\bar q}(t)\,dt
\end{align}
where $\dot {\bar q}(t)$ is the optimal intrinsic exercise strategy. It can be represented in the form
\begin{align}
    \label{eq:intrinsic_exercise}
    \dot {\bar q}(t) = r_{min} + \Delta r\,\theta[C-F(t)]
\end{align}
where 
\[
    \Delta r = r_{max} - r_{min} \qquad \quad
    \theta[x] = \left\{
        \begin{array}{ll}
            1 \qquad &  \text{for }x>0 \\
            0 \qquad &  \text{for }x<0
        \end{array}
    \right.
\]
We simplified the expression for the optimal exercise by excluding the boundary touch events.

Since the terminal level has to be preserved, the variation of both sides of Eq.~(\ref{eq:boundary_condition2}) must vanish. Substituting Eq.~(\ref{eq:intrinsic_exercise}) into Eq.~(\ref{eq:boundary_condition2}) and finding its first variation with respect to $C$ and $F$, we get
\begin{align}
    \int_0^{T_e} \Delta r \,\Big( \delta[C - F(t)]\;\delta C - 
        \delta[C - F(t)]\;\delta F(t) \Big)\,dt \approx 0
\end{align}
where $\delta[x] = \theta'[x]$ is the Dirac delta function. From here follows
\begin{align}
    \label{eq:dC_1_order}
    \delta C \approx \frac{\Delta r}{K} \int_0^{T_e} \delta[C - F(t)]\;\delta F(t)\,\,dt
    \qquad \text{where}\quad 
    K = \Delta r \int_0^{T_e} \delta[C - F(t)]\,\,dt
\end{align}

% \bibliographystyle{spmpsci}      % mathematics and physical sciences
\bibliographystyle{plain}
\bibliography{Storage_Option_Intrinsic}   % name your BibTeX data base

\end{document}